\documentclass[12pt]{amsart}

\usepackage{epsf}
\usepackage{verbatim}
\usepackage{amssymb}

\advance \topmargin by -\headheight
\advance \topmargin by -\headsep
     
\evensidemargin \oddsidemargin
\advance\hoffset by -3mm  
\advance\voffset by  8mm  

\theoremstyle{plain}
\newtheorem*{thm}{Theorem}
\newtheorem{lemma}{Lemma}

\theoremstyle{definition}
\newtheorem*{rk}{Remark}

\DeclareMathOperator{\End}{End}

\newcommand{\ra}{\rightarrow}

\newcommand{\C}{\mathbb C}

\newcommand{\Z}{\mathbb Z}

\newcommand{\skM}{\mbox{${\mathcal S}M$}}

\newcommand{\skH}{\mbox{${\mathcal S}H$}}
\newcommand{\skT}{\mbox{${\mathcal S}T$}}

\newcommand{\rM}{\mbox{${\mathcal R}M$}}

\newcommand{\rT}{\mbox{${\mathcal R}T$}}

\newcommand{\rH}{\mbox{${\mathcal R}H$}}

\begin{document}

\title{Irreducibility of some quantum representations of mapping class
groups} \author{Justin Roberts} \address{Department of Mathematics and
Statistics, Edinburgh University, EH3 9JZ, Scotland, UK}
\email{justin@maths.ed.ac.uk} \date{September 21, 1999.}

\begin{abstract}
The $SU(2)$ TQFT representation of the mapping class group of a closed
surface of genus $g$, at a root of unity of prime order, is shown to
be irreducible. Some examples of reducible representations are also
given.
\end{abstract}

\maketitle

\section{Introduction}

The Witten-Reshetikhin-Turaev topological quantum field theories (see
Reshetikhin and Turaev \cite{RT} or Turaev's book \cite{T}) provide
many interesting finite-dimensional representations of mapping class
groups of surfaces, about which little is currently known. In this
paper we will consider only the representations coming from the
$SU(2)$ theory, which may be defined and studied via the Kauffman
bracket skein theory (see Lickorish \cite{L}, Roberts
\cite{R1}). Calculations using this approach are easier and more
concrete than in the more general cases (where one has to work more
explicitly with quantum groups) but typically provide insight into the
general cases, which can be worked out along the same lines (as for
example with the integrality results of Masbaum-Roberts \cite{MR} and
then Masbaum-Wenzl \cite{MW}).

Recent papers by Funar \cite{Fun} and Masbaum \cite{M} studied the
question of whether the image of the mapping class group of a closed
surface under such a representation (at an $r$th root of unity) was
infinite or not. Here, another aspect will be considered: are the
representations irreducible? I was asked this question by Ivan Smith,
who was interested in the geometric quantization approach to the
representations, and found it difficult to answer using such an
algebro-geometric approach. Surprisingly, even the TQFT literature
does not seem to provide the answer. The purpose of this note,
therefore, is to make a start by answering the question at least in
some cases.

\begin{thm}
Let $r \geq 3$ be prime. Then the $SU(2)$ TQFT representation of the
mapping class group of a closed surface of genus $g$, at an $r$th root
of unity, is irreducible.
\end{thm}

Unfortunately the proof which will be explained below does not seem to
generalise in a straightforward way (unlike, say, the methods of
\cite{R2}) to either non-prime $r$ or to higher rank quantum groups. 
Neither is it completely clear how to extend it to the case of
surfaces with punctures.


\section{Skein theory preliminaries}

The proof of the theorem is not very long or complicated, so this
explanation of the background will be kept equally brief. The main
purpose is simply to fix the notation. The paper by Lickorish
\cite{L} and the book of Kauffman and Lins \cite{KL} contain full
explanations of how one uses skein theory to build $3$-manifold
invariants, whilst \cite{BHMV, R3} explain how to develop a full TQFT
using the same principles.

Fix $r \geq 3$, the integer `level', and let $A= e^{2\pi i/4r}$. The
symbol $\skM$ denotes the Kauffman skein space (a complex vector space
defined using the parameter $A$) of a compact oriented 3-manifold
$M$. It is the vector space generated by isotopy classes (rel
boundary) of framed links inside $M$, modulo the usual local Kauffman
bracket relations.

The Jones-Wenzl idempotents $f^{(a)}$ of the Temperley-Lieb algebras
(skein spaces of a cylinder with $a$ points at each end) are defined
for $a=0, 1, \ldots, r-1$. Inside any skein space $\skM$ one can
consider the subspace spanned by elements consisting of $f^{(r-1)}$ in
a small cylinder with the ends connected up in any way. Factoring out
by this subspace gives the {\em reduced skein space} $\rM$. It is a
fact (see Roberts \cite{R3}) that the reduced skein space of any
$3$-manifold depends only on its boundary, and gives a model for the
Witten-Reshetikhin-Turaev vector space of this boundary.

In particular if $H$ is a handlebody and $\Sigma$ its boundary then
the reduced skein space $\rH$ is identified with the W-R-T space
usually written (for example in Blanchet, Habegger, Masbaum and Vogel
\cite{BHMV}) as $V(\Sigma)$. In \cite{BHMV}, the space $V(\Sigma)$ is
constructed as a quotient of $\skH$, but in a fairly abstract way;
the point of the reduced skein space is simply that it is an explicit
local combinatorial description of the quotient.

The most important such spaces are those associated to a solid torus,
denoted $\skT$ and $\rT$ for convenience. The skein space $\skT$ is a
polynomial algebra with one generator $\alpha$, a single curve winding
once around the torus. The elements $\phi_a \in \skT$ ($a=0,1, \ldots,
r-2$), given by taking the $a$th Chebyshev polynomial of $\alpha$ (or
by closing up Jones-Wenzl idempotents) are particularly important, as
they descend to a basis for the $(r-1)$-dimensional space
$\rT$. Lickorish's construction of $3$-manifold invariants is based on
an element $\Omega \in \rT$ defined by
\[\Omega =\eta \sum_{a=0}^{r-2} \Delta(a) \phi_a,\]
where
\[ \eta= \frac{A^2 - A^{-2}}{i \sqrt{2r}} =\sqrt{{2/r}} \sin(\pi/r)\quad
\hbox{and}\ \Delta(a) = (-1)^a\frac{A^{2(a+1)}-A^{-2(a+1)}}{A^2 -
A^{-2}}.\]
(In \cite{BHMV} and \cite{MR}, the symbol $\omega$ was used instead of
$\Omega$; the convention here agrees with the one in \cite{R2}.)

For a handlebody $H$ of genus $g \geq 2$, one can again write down a
basis of $\rH$. The usual basis is given by picking $3g-3$ discs
chopping up $H$ in a pants decomposition, and then drawing a trivalent
spine dual to the decomposing discs. The standard basis elements $v$
are made by attaching Jones-Wenzl idempotents to the edges of this
graph and joining them suitably at the vertices. They are parametrised
by the labellings of their idempotents, in other words by a subset of
the set of labellings of the edges by integers in the range $0$ to
$r-2$. The vacuum vector $v_0$ is the basis vector corresponding to
the empty link (all labels are $0$).

The action of the mapping class group $\Gamma_g$ on $\rH$ can be
defined in a natural but implicit way (see \cite{R1}), but here it is
more useful to have an explicit description of the action of a
positive Dehn twist $T_{\gamma}$ about a curve $\gamma \in
\Sigma$. If $x$ is an element of $\rH$ then $T_\gamma x$ is
represented by adjoining to a skein element representing $x$ the curve
$\gamma$ with $\Omega$ inserted onto it (with framing $-1$ relative to
the surface), drawn just inside the boundary of $H$.

In particular, if $\gamma$ is the boundary of one of the pants discs,
the twist $T_\gamma$ acts on a standard basis vector $v$ by putting a
full positive twist in the edge passing through this disc. Idempotents
are eigenvectors under such twist operations. Consequently, if the
relevant edge is coloured with $a$, then $T_\gamma v = \xi_a v$, where
$\xi_a = (-1)^a A^{a^2+2a}$ is the associated eigenvalue (twist
coefficient).

\begin{rk}
The representation defined by these twist generators (or as
constructed in \cite{R1}) is only projective, and it is customary to
lift to a genuine linear representation of a central extension of
$\Gamma_g$. However, the projective ambiguity has no bearing on the
question of irreducibility, so this fact can be safely ignored.
\end{rk}


\section{Proof of the theorem}

\begin{lemma}
If $r$ is prime then the vectors $t_b$, for $b=0, 1, \ldots, r-2$,
given by placing $b$ parallel $-1$-framed copies of $\Omega$ in the
solid torus form a basis for $\rT$.
\end{lemma}
\begin{proof}
The quickest way to see this is to use the non-degenerate pairing $\rT
\times \rT \ra \C$ obtained by gluing together two solid tori to make
$S^3$ (whose skein space is canonically $\C$). Pairing $t_b$ with
$\phi_a$ results in $\xi_b^a \Delta(a)$, so that the change of basis
matrix expressing the $t_b$, viewed as linear functionals, in terms of
the dual functionals $\phi_a^*$ is a Vandermonde matrix $(\xi_a^b)$
times a diagonal matrix whose diagonal entries are the (non-zero)
$\Delta(a)$.  Now $\xi_a= (-1)^a A^{a^2+2a}$, and one can easily check
that these are all distinct when $r$ is prime, hence the first matrix
is invertible. The second is obviously invertible, and so the $t_b$
indeed form a basis.
\end{proof}

\begin{lemma}
Suppose $C$ is a collection of disjoint curves on $\Sigma_g$. Then one
can obtain an element of $\rH_g$ by viewing these curves as lying in
$H_g$ and attaching $\Omega$ to each one with framing $-1$ relative to
the surface. Such vectors $\Omega(C)$ span $\rH_g$. Therefore the
image of the group algebra of $\Gamma_g$ applied to the vacuum vector
$v_0$ is all of $\rH_g$.
\end{lemma}
\begin{proof}
Consider $H_g$ as a thickened $(g-1)$-holed disc.  The reduced skein
space is spanned by a finite number of elements of $\skH_g$, which may
be represented by links lying in the holed disc (consider the usual
planar projection onto this disc). Given such a planar link $L$,
isotop it to be near to the boundary $\Sigma_g$, and rewrite each of
its components (thought of as an element of $\rT$) as a linear
combination of the basis elements $t_b$. Since each $t_b$ is really
$\Omega(C)$, for $C$ a single curve parallelled $b$ times, this
proves the assertion about spanning. The final part follows
immediately from the description of the action of a Dehn twist on
$\rH_g$.
\end{proof}

Note that this does not immediately imply irreducibility. For example,
if $\Z$ acts on $\C^2$ with distinct eigenvalues then the action is
reducible but the group algebra applied to any non-eigenvector is all
of $\C^2$. (This example can even be chosen to be unitary.)

\begin{lemma}
The subgroup $P_g$ of the mapping class group generated by Dehn twists
on the standard pants curves is a free abelian group, under whose
action $\rH_g$ breaks up as a sum of one-dimensional representations,
spanned by the standard basis vectors $v$.
\end{lemma}
\begin{proof}
The standard basis vectors are certainly simultaneous eigenvectors for
the twists generating $P_g$, as each twist just multiplies the vector
by a twist coefficient $\xi_a$. But their collections of eigenvalues
are distinct since the $\xi_a$ are (when $r$ is prime) and hence they
span individual one-dimensional eigenspaces.
\end{proof}

Now the theorem can be proved. Suppose $\theta$ is any endomorphism of
$\rH_g$ commuting with the action of the whole mapping class group
$\Gamma_g$. Then it certainly acts diagonally with respect to the
standard basis, because it commutes in particular with the subgroup
$P_g$ whose eigenvectors they are. Let us write $\theta v = \lambda_v
v$ for a standard basis vector $v$. All we need to do is to show that
the matrix of $\theta$ is actually a {\em scalar} to conclude, via
Schur's lemma, that $\rH_g$ is an irreducible representation of
$\Gamma_g$. To see this, observe that any standard basis vector $v$
can be generated from the vacuum vector $v_0$ by the action $\psi$ of
some element of the group algebra of $\Gamma_g$, by lemma 2. Then,
since $\theta$ commutes with $\psi$,
\[ \theta v= \theta \psi v_0 = \psi \theta v_0 = \lambda_{v_0} \psi v_0 =
\lambda_{v_0} v,\]
but also $\theta v = \lambda_v v$, so $\theta$ is a scalar. 

\section{Further comments}

There are certainly cases in which the representations are not
irreducible. It is difficult to find these in general, but in genus
$1$, there is a large body of literature studying modular invariant
partition functions for affine Lie algebras which provides a more than
complete solution. See Capelli-Itzykson-Zuber \cite{CIZ}, Fuchs' book
\cite{Fuc}, and Gannon \cite{G} for a short proof of the result of
\cite{CIZ}.

The problem studied in these references is to find all
$SL(2,\Z)$-invariant linear combinations
\[ Z = \sum_{a,b=0}^{r-2} Z_{a,b} \phi_a \otimes \bar\phi_b \in V
\otimes \bar V,\]where $V=V(\Sigma_1)$ (and $SL(2,\Z)=\Gamma_1$ is the
mapping class group of the torus acting on it), and the coefficients
$Z_{a,b}$ are non-negative integers. Since $\bar V \cong V^*$, such
elements can be thought of as invariant endomorphisms of $V$, and a
reasonable first step in classifying such elements $Z$ is to find the
commutant of $SL(2,\Z)$ in $\End(V)$. This is carried out by Cappelli,
Itzykson and Zuber \cite{CIZ}, who find the dimension of the commutant
in terms of the number of divisors of $r$. (They then go on to find
the non-negative integer matrices $Z$ lying in the commutant and to
obtain an amazing $A-D-E$ classification.) In particular, the
commutant is trivial when $r$ is prime, which agrees with the
irreducibility theorem proved above, and also shows that it is sharp
in genus $1$. 

The higher-genus situation does not seem to have been studied much. It
is worth observing that the method of \cite{CIZ} for finding the
commutant in genus $1$ relies on averaging over the image of
$\Gamma_1=SL(2,Z)$, which is finite (see for example Gilmer
\cite{Gi}). In higher genus, as noted in the introduction, the image is
known to be (with finitely many small $r$ exceptions) infinite, so
such methods fail. Whether one can apply the methods of skein theory
to the problem of modular invariants, or vice versa, remains to be
seen.


\end{document}